\documentclass[12pt]{amsart}
\usepackage{amssymb}
\begin{document}
\theoremstyle{plain}
\newtheorem{thm}{Theorem}[section]
\newtheorem{lem}[thm]{Lemma}
\newtheorem{cor}[thm]{Corollary}
\newtheorem{prop}[thm]{Proposition}

\theoremstyle{definition}
\newtheorem{rem}[thm]{Remark}
\newtheorem{defn}[thm]{Definition}
\newtheorem{ex}[thm]{Example}

\numberwithin{equation}{section}
\newcommand{\Hom}{{\rm Hom}}
\newcommand{\End}{{\rm End}}
\newcommand{\res}{{\rm res}}
\newcommand{\sA}{{\mathcal A}}
\newcommand{\sB}{{\mathcal B}}
\newcommand{\sC}{{\mathcal C}}
\newcommand{\sD}{{\mathcal D}}
\newcommand{\sE}{{\mathcal E}}
\newcommand{\sF}{{\mathcal F}}
\newcommand{\sG}{{\mathcal G}}
\newcommand{\sH}{{\mathcal H}}
\newcommand{\sI}{{\mathcal I}}
\newcommand{\sJ}{{\mathcal J}}
\newcommand{\sK}{{\mathcal K}}
\newcommand{\sL}{{\mathcal L}}
\newcommand{\sM}{{\mathcal M}}
\newcommand{\sN}{{\mathcal N}}
\newcommand{\sO}{{\mathcal O}}
\newcommand{\sP}{{\mathcal P}}
\newcommand{\sQ}{{\mathcal Q}}
\newcommand{\sR}{{\mathcal R}}
\newcommand{\sS}{{\mathcal S}}
\newcommand{\sT}{{\mathcal T}}
\newcommand{\sU}{{\mathcal U}}
\newcommand{\sV}{{\mathcal V}}
\newcommand{\sW}{{\mathcal W}}
\newcommand{\sX}{{\mathcal X}}
\newcommand{\sY}{{\mathcal Y}}
\newcommand{\sZ}{{\mathcal Z}}
\newcommand{\A}{{\mathbb A}}
\newcommand{\B}{{\mathbb B}}
\newcommand{\C}{{\mathbb C}}
\newcommand{\D}{{\mathbb D}}
\newcommand{\E}{{\mathbb E}}
\newcommand{\F}{{\mathbb F}}
\newcommand{\G}{{\mathbb G}}
\renewcommand{\H}{{\mathbb H}}
\newcommand{\I}{{\mathbb I}}
\newcommand{\J}{{\mathbb J}}
\newcommand{\M}{{\mathbb M}}
\newcommand{\N}{{\mathbb N}}
\renewcommand{\P}{{\mathbb P}}
\newcommand{\Q}{{\mathbb Q}}
\newcommand{\R}{{\mathbb R}}
\newcommand{\T}{{\mathbb T}}
\newcommand{\U}{{\mathbb U}}
\newcommand{\V}{{\mathbb V}}
\newcommand{\W}{{\mathbb W}}
\newcommand{\X}{{\mathbb X}}
\newcommand{\Y}{{\mathbb Y}}
\newcommand{\Z}{{\mathbb Z}}


\catcode`\@=11
%
%
\def\opn#1#2{\def#1{\mathop{\kern0pt\fam0#2}\nolimits}} 
\def\bold#1{{\bf #1}}%
\def\underrightarrow{\mathpalette\underrightarrow@}
\def\underrightarrow@#1#2{\vtop{\ialign{$##$\cr
 \hfil#1#2\hfil\cr\noalign{\nointerlineskip}%
 #1{-}\mkern-6mu\cleaders\hbox{$#1\mkern-2mu{-}\mkern-2mu$}\hfill
 \mkern-6mu{\to}\cr}}}
\let\underarrow\underrightarrow
\def\underleftarrow{\mathpalette\underleftarrow@}
\def\underleftarrow@#1#2{\vtop{\ialign{$##$\cr
 \hfil#1#2\hfil\cr\noalign{\nointerlineskip}#1{\leftarrow}\mkern-6mu
 \cleaders\hbox{$#1\mkern-2mu{-}\mkern-2mu$}\hfill
 \mkern-6mu{-}\cr}}}
%
%

%
\def\:{\colon}
\let\oldtilde=\tilde
\def\tilde#1{\mathchoice{\widetilde{#1}}{\widetilde{#1}}%
{\indextil{#1}}{\oldtilde{#1}}}
\def\indextil#1{\lower2pt\hbox{$\textstyle{\oldtilde{\raise2pt%
\hbox{$\scriptstyle{#1}$}}}$}}
\def\pnt{{\raise1.1pt\hbox{$\textstyle.$}}}
%

%
\let\amp@rs@nd@\relax
\newdimen\ex@
\ex@.2326ex
\newdimen\bigaw@
\newdimen\minaw@
\minaw@16.08739\ex@
\newdimen\minCDaw@
\minCDaw@2.5pc
\newif\ifCD@
\def\minCDarrowwidth#1{\minCDaw@#1}
\newenvironment{CD}{\@CD}{\@endCD}
\def\@CD{\def\A##1A##2A{\llap{$\vcenter{\hbox
 {$\scriptstyle##1$}}$}\Big\uparrow\rlap{$\vcenter{\hbox{%
$\scriptstyle##2$}}$}&&}%
\def\V##1V##2V{\llap{$\vcenter{\hbox
 {$\scriptstyle##1$}}$}\Big\downarrow\rlap{$\vcenter{\hbox{%
$\scriptstyle##2$}}$}&&}%
\def\={&\hskip.5em\mathrel
 {\vbox{\hrule width\minCDaw@\vskip3\ex@\hrule width
 \minCDaw@}}\hskip.5em&}%
\def\verteq{\Big\Vert&&}%
\def\noarr{&&}%
\def\vspace##1{\noalign{\vskip##1\relax}}\relax\let\amp@rs@nd@&\iffalse}\fi
 \CD@true\vcenter\bgroup\relax\let\\=\cr\iffalse}\fi\tabskip\z@skip\baselineskip20\ex@
 \lineskip3\ex@\lineskiplimit3\ex@\halign\bgroup
 &\hfill$\m@th##$\hfill\cr}
\def\@endCD{\cr\egroup\egroup}
%
\def\>#1>#2>{\amp@rs@nd@\setbox\z@\hbox{$\scriptstyle
 \;{#1}\;\;$}\setbox\@ne\hbox{$\scriptstyle\;{#2}\;\;$}\setbox\tw@
 \hbox{$#2$}\ifCD@
 \global\bigaw@\minCDaw@\else\global\bigaw@\minaw@\fi
 \ifdim\wd\z@>\bigaw@\global\bigaw@\wd\z@\fi
 \ifdim\wd\@ne>\bigaw@\global\bigaw@\wd\@ne\fi
 \ifCD@\hskip.5em\fi
 \ifdim\wd\tw@>\z@
 \mathrel{\mathop{\hbox to\bigaw@{\rightarrowfill}}\limits^{#1}_{#2}}\else
 \mathrel{\mathop{\hbox to\bigaw@{\rightarrowfill}}\limits^{#1}}\fi
 \ifCD@\hskip.5em\fi\amp@rs@nd@}
\def\<#1<#2<{\amp@rs@nd@\setbox\z@\hbox{$\scriptstyle
 \;\;{#1}\;$}\setbox\@ne\hbox{$\scriptstyle\;\;{#2}\;$}\setbox\tw@
 \hbox{$#2$}\ifCD@
 \global\bigaw@\minCDaw@\else\global\bigaw@\minaw@\fi
 \ifdim\wd\z@>\bigaw@\global\bigaw@\wd\z@\fi
 \ifdim\wd\@ne>\bigaw@\global\bigaw@\wd\@ne\fi
 \ifCD@\hskip.5em\fi
 \ifdim\wd\tw@>\z@
 \mathrel{\mathop{\hbox to\bigaw@{\leftarrowfill}}\limits^{#1}_{#2}}\else
 \mathrel{\mathop{\hbox to\bigaw@{\leftarrowfill}}\limits^{#1}}\fi
 \ifCD@\hskip.5em\fi\amp@rs@nd@}
%
%
\newenvironment{CDS}{\@CDS}{\@endCDS}
\def\@CDS{\def\A##1A##2A{\llap{$\vcenter{\hbox
 {$\scriptstyle##1$}}$}\Big\uparrow\rlap{$\vcenter{\hbox{%
$\scriptstyle##2$}}$}&}%
\def\V##1V##2V{\llap{$\vcenter{\hbox
 {$\scriptstyle##1$}}$}\Big\downarrow\rlap{$\vcenter{\hbox{%
$\scriptstyle##2$}}$}&}%
\def\={&\hskip.5em\mathrel
 {\vbox{\hrule width\minCDaw@\vskip3\ex@\hrule width
 \minCDaw@}}\hskip.5em&}
\def\verteq{\Big\Vert&}
\def\novarr{&}
\def\noharr{&&}
\def\SE##1E##2E{\slantedarrow(0,18)(4,-3){##1}{##2}&}
\def\SW##1W##2W{\slantedarrow(24,18)(-4,-3){##1}{##2}&}
\def\NE##1E##2E{\slantedarrow(0,0)(4,3){##1}{##2}&}
\def\NW##1W##2W{\slantedarrow(24,0)(-4,3){##1}{##2}&}
\def\slantedarrow(##1)(##2)##3##4{%
\thinlines\unitlength1pt\lower 6.5pt\hbox{\begin{picture}(24,18)%
\put(##1){\vector(##2){24}}%
\put(0,8){$\scriptstyle##3$}%
\put(20,8){$\scriptstyle##4$}%
\end{picture}}}
\def\vspace##1{\noalign{\vskip##1\relax}}\relax\let\amp@rs@nd@&\iffalse}\fi
 \CD@true\vcenter\bgroup\relax\let\\=\cr\iffalse}\fi\tabskip\z@skip\baselineskip20\ex@
 \lineskip3\ex@\lineskiplimit3\ex@\halign\bgroup
 &\hfill$\m@th##$\hfill\cr}
\def\@endCDS{\cr\egroup\egroup}
%
\newdimen\TriCDarrw@
\newif\ifTriV@
\newenvironment{TriCDV}{\@TriCDV}{\@endTriCD}
\newenvironment{TriCDA}{\@TriCDA}{\@endTriCD}
\def\@TriCDV{\TriV@true\def\TriCDpos@{6}\@TriCD}
\def\@TriCDA{\TriV@false\def\TriCDpos@{10}\@TriCD}
\def\@TriCD#1#2#3#4#5#6{%
\setbox0\hbox{$\ifTriV@#6\else#1\fi$}
\TriCDarrw@=\wd0 \advance\TriCDarrw@ 24pt
\advance\TriCDarrw@ -1em
\def\SE##1E##2E{\slantedarrow(0,18)(2,-3){##1}{##2}&}
\def\SW##1W##2W{\slantedarrow(12,18)(-2,-3){##1}{##2}&}
\def\NE##1E##2E{\slantedarrow(0,0)(2,3){##1}{##2}&}
\def\NW##1W##2W{\slantedarrow(12,0)(-2,3){##1}{##2}&}
\def\slantedarrow(##1)(##2)##3##4{\thinlines\unitlength1pt
\lower 6.5pt\hbox{\begin{picture}(12,18)%
\put(##1){\vector(##2){12}}%
\put(-4,\TriCDpos@){$\scriptstyle##3$}%
\put(12,\TriCDpos@){$\scriptstyle##4$}%
\end{picture}}}
\def\={\mathrel {\vbox{\hrule
   width\TriCDarrw@\vskip3\ex@\hrule width
   \TriCDarrw@}}}
\def\>##1>>{\setbox\z@\hbox{$\scriptstyle
 \;{##1}\;\;$}\global\bigaw@\TriCDarrw@
 \ifdim\wd\z@>\bigaw@\global\bigaw@\wd\z@\fi
 \hskip.5em
 \mathrel{\mathop{\hbox to \TriCDarrw@
{\rightarrowfill}}\limits^{##1}}
 \hskip.5em}
\def\<##1<<{\setbox\z@\hbox{$\scriptstyle
 \;{##1}\;\;$}\global\bigaw@\TriCDarrw@
 \ifdim\wd\z@>\bigaw@\global\bigaw@\wd\z@\fi
 \mathrel{\mathop{\hbox to\bigaw@{\leftarrowfill}}\limits^{##1}}
 }
 \CD@true\vcenter\bgroup\relax\let\\=\cr\iffalse}\fi
 \tabskip\z@skip\baselineskip20\ex@
 \lineskip3\ex@\lineskiplimit3\ex@
 \ifTriV@
 \halign\bgroup
 &\hfill$\m@th##$\hfill\cr
#1&\multispan3\hfill$#2$\hfill&#3\\
&#4&#5\\
&&#6\cr\egroup%
\else
 \halign\bgroup
 &\hfill$\m@th##$\hfill\cr
&&#1\\%
&#2&#3\\
#4&\multispan3\hfill$#5$\hfill&#6\cr\egroup
\fi}
\def\@endTriCD{\egroup}

\title[Semistable bundles and irreducibility]
{Semistable bundles on curves and irreductible representations
of the fundamental group} 
\author{ H\'el\`ene Esnault }
\author{ Eckart Viehweg} 
\address{ Universit\"at Essen, FB6 Mathematik, 45 117 Essen, Germany}
\email{ esnault@uni-essen.de}
\email{ viehweg@uni-essen.de}
\thanks{ This work has been partly supported by the DFG Forschergruppe
''Arithmetik und Geometrie''}

\maketitle
\setcounter{section}{-1} 

\section{Introduction}

In this note, we make an attempt to understand the meaning of
Bolibruch's theorem for curves of higher genus. 

\begin{thm}[Bolibruch \cite{Bo}] \label{thm:bo} 
Let $$\rho : \pi_1 (\P^{1}_{\C} - \Sigma) \>>> GL(N,\C)$$
be an irreducible representation of the fundamental group of the
complement of finitely many points $\Sigma \neq \emptyset$. Then
there is a logarithmic connection 
$$
\nabla :  \sO^{\oplus{N}} \>>> \Omega_X^1 ({\rm log} \ \Sigma)
\otimes \sO^{\oplus{N}}
$$
such that the local system $\ker(\nabla|_{X-\Sigma})$ on
$\P^{1}_{\C} - \Sigma$ is defined by $\rho$. 
\end{thm}
Bolibruch's proof is very analytic, but Gabber (\cite{G}) gave a
more algebraic approach, which we recall in section
\ref{sec:gab} (see also \cite{Le}). Using his construction, we interpret 
Bolibruch's theorem in the following way.
\begin{thm} \label{thm:ev} 
Let $X$ be a curve over an algebraically closed field $k$ of
characteristic 0, and let $\emptyset \neq \Sigma \subset X (k)$ 
consist of finitely many points. Let 
$$
\nabla : E \>>> \Omega_X^1 ({\rm log} \ \Sigma ) \otimes E
$$
be a logarithmic connection on a vectorbundle $E$ of rank $N$
such that for all subsheaves $\{ 0\} \neq F
\subset E$ with ${\rm rank}(F) < N$,
$$
\nabla F \not\subset \Omega_X^1 ({\rm log} \ \Sigma ) \otimes F. 
$$
Then for any $p \in \Sigma$, there is a semistable vectorbundle
$E'$ of degree $0$ and a logarithmic connection 
$$
\nabla' : E' \>>> \Omega_X^1 ({\rm log} \ \Sigma) \otimes E', 
$$
with $(E',\nabla')|_{X-\{p\}} = (E,\nabla)|_{X-\{p\}}$. 
\end{thm}
Any semistable bundle $E'$ of rank $N$ and degree 0, has a
canonical filtration (see (\ref{socle})), the
graded bundles $gr_iE'$ of which are direct sums of stable ones. 
Due to the Narasimhan-Seshadri correspondence \cite{NS} 
over $\C$, there is a unitary connection $d_i$
on $grE'_i$ which is uniquely defined. 

The curious point is that, over $k=\C$, we associate 
to an irreducible representation of the fundamental group 
$$
\pi_1(X-\Sigma) \>>> GL(N, \C)
$$
of the open curve $X-\Sigma$, unitary representations of the 
fundamental group of the compact curve
$$
\pi_1(X)\>>> U(N_i,
\C),\mbox{ \ where \ } \sum_iN_i=N,
$$ 
via theorem \ref{thm:ev} and the Narasimhan-Seshadri
correspondence.
 
Conversely it is easy to associated such unitary representations
of $\pi_1(X)$ an irreducible representation $\pi_1(X-\Sigma) \to
GL(N, \C)$: 
\begin{prop} \label{thm:ev2} Let $X$ be a curve over $\C$ 
let $E$ be a semistable bundle on $X$ of degree 0 with graded
bundles $gr_i(E)$ for the canonical filtration. 
\begin{enumerate}
\item[1)] There is a connection $\nabla: E \to \Omega^1_X
\otimes E$ respecting the canonical filtration  on
$E$, such that $gr_i(\nabla) = d_i$.
\item[2)] There is a constant $\sigma \leq 3$ depending only on $E$
such that for any reduced divisor $\Sigma$ with ${\rm deg}( 
\Sigma) \geq \sigma$, there is a connection 
$$\nabla: E \>>> \Omega^1_X(\log \Sigma) \otimes E$$
such that for all subsheaves $\{0\} \neq F
\subset E$ with ${\rm rank}(F) < N$,
$$
\nabla F \not\subset \Omega_X^1 ({\rm log} \ \Sigma ) \otimes F. 
$$
\end{enumerate}
\end{prop}
This  way of going back and forth between representations of the
projective and the open curve is very lose. On both
sides one has parameters. It is not clear whether one should
think of this really as a correspondence. It is also not clear
how to interpret this in terms of compactification of the moduli
space of stable bundles of degree 0. 

\section{Gabber's construction }  \label{sec:gab}

We explain Gabber's construction, transposing it to the
algebraic context of theorem \ref{thm:ev}. Hence we consider
a projective curve $X$ over $k$, a divisor
$\Sigma >0$ and a logarithmic connection
$$
\nabla : E \>>> \Omega_X^1 ({\rm log} \ \Sigma) \otimes E
$$
on a vectorbundle $E$. We fix a point $p \in \Sigma$ and
denote by 
$$
\Gamma=\res_p(\nabla):E\otimes k(p) \>>> E\otimes k(p)
$$
the residue of $\nabla$. 

For $0 \neq w \in E \otimes k (p)$ define $E'_w$ to be the
inverse image of $k_w$ under the restriction map 
$E \to E \otimes k (p),$
and $E_w = E'_w (p)$. Then $E \subset E_w \subset E (p)$ and ${\rm
deg} \ E_w = {\rm deg}(E)+1$. 

The connection $\nabla $ extends to
$\nabla_w$ on $E_w$ if and only if $w$ is an eigenvector of
$\Gamma$. More precisely, let $(w, e_2, \ldots , e_N)$ be a
basis of $E \otimes k (p)$ in which $\Gamma = (\gamma_{ij})$ is
triangular, that is $\gamma_{ij} =0 \ i >j$.
Then in the basis $(\frac{w}{t}, e_2,
\ldots , e_N)$ of $E_w \otimes k (p)$ the residue
$\res_p(\nabla_w)=\Gamma_w = (\gamma'_{ij})$ fulfills: 
\begin{align*}
\gamma'_{ij} &= \gamma_{ij} \mbox{ \ for \ } i \geq 2, j \geq 2\\
\gamma'_{11} &= \gamma_{11} - 1 \\
\gamma'_{1i} &=0 \ i \geq 2 .
\end{align*}
Thus the roots of the characteristic polynomial of $\Gamma_w$,
are $\gamma_{11} -1 , \gamma_{22} , \ldots , \gamma_{NN}$. 

\begin{defn} We say that $(E', \nabla')$ is obtained from $(E,
\nabla)$ by an elementary $G$-transformation at $p$ if there is
an eigenvector $0 \neq w \in E \otimes k (p)$ of $\Gamma$ such
that $(E', \nabla') = (E_w, \nabla_w)$. 
\end{defn}

\begin{thm}[Gabber] \label{gabber}
Let $\nabla : E \to \Omega_X^1 ({\rm log} \
\Sigma) \otimes E$ be any connection, and $M \in \N$.
Then there is a connection 
$$
\nabla' : E' \>>> \Omega_X^1 ({\rm log} \ \Sigma) \otimes E'
$$
with $(E' , \nabla') |_{X - \{ p \}} = (E, \nabla) |_{X - \{ p \}}$
such that 
\begin{enumerate}
\item[1)] the characteristic polynomial of
$\Gamma'=\res_p(\nabla')$ has no multiple zeros,
\item[2)] if $\lambda, \mu$ are 2 eigenvalues of $\Gamma'$, with
$\lambda - \mu \in \Z$, then $|\lambda - \mu| \geq M$,  
\item[3)] $(E' , \nabla')$ is obtained from $(E, \nabla)$ by at
most $\frac{N^3 M}{2}$ elementary $G$-trans\-for\-ma\-tions at $p$. 
\end{enumerate}
\end{thm} 
\begin{proof} One orders the roots of the characteristic
polynomial of $\Gamma$ in subsets
$I_1, \ldots , I_\ell$, 
$$
I_j = \{ \lambda_{j,1}, \ldots , \lambda_{j,mj} \} ,
\mbox{ \ where \ }\sum^{\ell}_{j=1} m_j = N
$$
such that $0 \leq \lambda_{j, i+1} - \lambda_{j,i} \in \N$, and
$
\lambda_{j,s} - \lambda_{j',s'} \not\in \Z \ \mbox{for} \  j'
\neq j .
$
By taking an eigenvector $e_1\in E\otimes k(p)$ for
$\lambda_{11}$ and replacing $E$ by $E_{e_1}$, one transforms $I_1$ to 
$$
I_1 = \{ \lambda_{1,1} -1, \lambda_{1,2}, \ldots , \lambda_{1,
m_1} \} .
$$
Repeating this $m_1 M$ times, one replaces $I_1$ by
$$
I_1 = \{ \lambda_{1,1} - m_1 M , \lambda_{1,2} , \ldots ,
\lambda_{1, m_1} \}.
$$
Since $\lambda_{1,1} - m_1 M\neq \lambda_{1,2}$, there exists
an eigenvector $e_2$ with eigenvalue $\lambda_{1,2}$, and repeating
the same transformation $(m_1 - 1) M$ times with $e_2$ instead
of $e_1$ one transforms $\lambda_{1,2}$ to
$\lambda_{1,2} - (m_1 -1)M$, without changing the other roots of
the characteristic polynomial.
After $\frac{m_1(m_1 -1)}{2} M$ steps, one has 
$$
I_1 = \{ \lambda_{1,1} - m_1 M , \lambda_{1,2} - (m_1 -1) M ,
\ldots , \lambda_{1, m_1-1}-M, \lambda_{1, m_1} \} .
$$
Repeating this for $I_2, \ldots , I_\ell$, one needs at most 
$$
(\sum^{\ell}_{j=1} \frac{m_j (m_j -1)}{2}) M \leq \frac{N^3}{2}
\cdot M $$
steps to satisfy the first and second condition in \ref{gabber}. 
\end{proof} 

\section{The proof of theorem \ref{thm:ev}}

Let $E_0 = 0 \subsetneqq E_1 \subset E_2 \subset \ldots
\subset E_m = E$ be the Harder-Narasimhan filtration \cite{HN}
of a rank $N$ vector bundle $E$,
uniquely determined by the two conditions: 
$$
\mu_i = \mu (E_i /E_{i-1}) < \mu_{i-1}
$$
and $E_i / E_{i-1}$ semistable,
where $\mu (F) = {\rm deg} (F) / {\rm rank} (F)$ for any vector
bundle.

In order to prove theorem \ref{thm:ev} we are allowed to
replace $E$ by $E(\ell p)$ for $\ell \in \Z$.
In fact, $\nabla$ stabilizes $E(\ell p)$ and the residue
$\Gamma$ of $\nabla$ in $p$ is replaced by $\Gamma - \ell {\rm Id}$. 
In particular this does not change the difference between
two eigenvalues of $\Gamma$. Thus, replacing $E$ by $E(\ell p)$,
we may assume that $- 1 < \mu (E_1) \leq 0$ and
consequently that ${\rm deg}(E) \leq 0 .$

\begin{lem} \label{spacing} If $\nabla : E \to \Omega_X^1 ({\rm log} \ \Sigma )
\otimes E$ does not stabilize any subbundle, and $-1 < \mu (E_1)
\leq 0$, then 
$$
- N - N^2 (2g-2+\sigma) \leq {\rm deg}(E) \leq 0
$$
where $g =$ genus of $X$, and $\sigma = |\Sigma |$. 
\end{lem} 
\begin{proof}
Let $i_0$ to be the minimal $i$ such that the
map 
$$
\eta_0: E_i \>>> \Omega_X^1 ({\rm log} \ \Sigma) \otimes E / E_{m-1}
$$
is not 0. Since $\nabla$ does not stabilize any subbundle, $i_0 \leq
m-1$, thus $\eta_0$ is linear and factors through $E_{i_0} /
E_{i_0-1}$. This shows that $\mu_{i_0} \leq \mu_{m} + (2g -2 +
\sigma)$. By assumption $\nabla$ does not 
stabilize $E_{i_0-1}$. Hence there exists some minimal number 
$i_1 \leq i_0-1$ such that
$$
\eta_1:E_i \>>> \Omega_X^1 ({\rm log} \ \Sigma) \otimes E / E_{i_0-1}
$$
is not trivial. Then $\eta_1$ factors through a linear map
$$ E_{i_1}/E_{i_1-1} \>>> \Omega_X^1 ({\rm log} \ \Sigma) \otimes 
E_{j_1} / E_{j_1-1}$$
for some $j_1$ with $i_0 \leq j_1 \leq m$. Consequently
$$
\mu_{i_1}\leq \mu_{j_1} + (2g-2+\sigma) 
\leq \mu_{i_0} + (2g-2+\sigma)\leq \mu_{m} + 2(2g-2+\sigma)  .
$$ 
One obtains inductively
$$
-1 \leq \mu_1 \leq \mu_m + m(2g-2+\sigma),
$$
and, since $\mu(E) \geq \mu_m$ and $N \geq m$, the inequality of lemma
\ref{spacing}. 
\end{proof}
Finally, one proves theorem \ref{thm:ev} in the following more
precise form:
\begin{thm} Let $(X, E, \nabla, \Sigma)$ be as in
theorem \ref{thm:ev}. Assume that 
$$
-1 < \mu (E_1) \leq 0,
$$
that the characteristic polynomial of 
$\Gamma=\res_p(\nabla)$ has no multiple zeros, and 
that 
$$
|\lambda - \mu | \geq M = N + N^2 (2g -2+\sigma)
$$
for different eigenvalues $\lambda$ and $\mu$
of $\Gamma$ with $\lambda - \mu \in
\Z$. 

Then there is a semistable vector bundle $E'$ of degree 0,
and an extension $\nabla'$ of $\nabla$ to $E'$, such that $(E',
\nabla')$ is obtained from $(E, \nabla)$ by at most $M$
elementary $G$-transformations at $p$. 
\end{thm} 

\begin{proof} 
We argue by induction on $- {\rm deg}(E)$ which is smaller than 
or equal to $M$
by lemma \ref{spacing}.

If ${\rm deg}(E) = 0$, $\mu(E_1)= \mu(E)=0$ as $\mu(E_1) \geq
\mu(E)$. Thus $E_1=E$ and $E$ is semistable of degree 0. 

Assume now that ${\deg}(E) < 0$. 
If $\mu(E_1) < 0$ as well, then for any elementary $G$
transformation at $p$, and any subsheaf $M\subset E_w$, one has 
$$
{\deg}(M) \leq {\deg} (M\cap E) + 1 \leq 0,
$$
thus 
$$ {\rm deg}(E_w)= {\rm deg}(E) +1 
\mbox{ \ and \ } -1 \le \mu((E_w)_1) \leq 0 .$$

Otherwise, $\mu(E_1)= {\rm deg}(E_1) =0$. We set $F=E_1$ for
notational simplicity and denote by $Q$ the
quotient $Q=E/F$. We consider an elementary $G$ transformation
at $p$ such that the eigenvector $w \in E\otimes k(p)$ maps non-trivially
to $Q\otimes k(p)$. One obtains an exact sequence
$$0\>>> F\>>> E_w \>>> Q_w \>>> 0 .$$
Let $(E_w)_1$ be the first bundle in the Harder-Narasimhan
filtration of $E_w$.
One certainly has 
$$
-1  \le \mu(E_1)\le \mu((E_w)_1).
$$
The inequality $\mu((E_w)_1) \leq 0$ is equivalent to the property that
${\rm deg}(M) \leq 0$ for all subsheaves $M \subset E_w$.
Consider $M \subset E_w$ 
and $M \subset M' \subset E_w$, where $M'$ is the inverse
image of $M/F\cap M$ under the projection $E_w \to E_w/F\cap M$. As 
$F \cap M \subset F$, one has ${\rm deg }(F \cap M) \leq 0$.
Thus
$${\rm deg}(M) \leq {\rm deg}(M/F \cap M) = {\rm deg}(M') + {\rm
deg}(F)={\rm deg}(M').$$
By definition of $E_1=F$, one has
$\mu ((E/F)_1) < 0$
and 
$$
{\rm deg}((M/ F\cap M)\cap Q)  \leq -1.
$$
This shows that
$${\rm deg}(M) \leq {\rm deg} (M/ F\cap M) \leq {\rm deg} \
((M/ F\cap M)\cap Q) +1 \leq -1 +1 \leq 0.$$
Thus again 
$$ {\rm deg}(E_w) = {\rm deg}(E) +1 \mbox{ \ and \ }
-1 < \mu((E_w)_1) \leq 0. $$
By induction we obtain the theorem.
\end{proof} 

\section{Existence of connections}

In this section we lift the unitary connections of the graded 
pieces of the canonical filtration.

\begin{lem}[Compare with \cite{Si}, lemma 3.5] \label{obs}
Let $X$ be an algebraic variety over a field $k$, 
$$0 \>>> S \> \iota >> E \> p >> Q \>>> 0 $$
be an extension of vector bundles
given by $u \in H^1 (X, \sH om (Q, S))$. Let $d_S$ and $d_Q$ be 
connections on $S$ and $Q$, respectively. 

Then there exists a connection
$\nabla $ on $E$ lifting $d_S$ and $d_Q$ if and only if $0 = du
\in H^1 (X, \Omega^{1}_{X} \otimes \sH om (Q, S))$, where $d =
\sH om (d_Q , d_S)$. Two such connections differ by an element
in $H^0 (X, \Omega^1_X \otimes \sH om (Q, S))$. 

In particular, if $X$ is projective smooth,
$k = \C$, and if $d$ is unitary, then $\nabla$ exists.  
\end{lem}

\begin{proof} Let $X = \bigcup U_i $ be an affine covering of $X$,
$$
\sigma_i : Q |_{U_i} \>>> E |_{U_i} , \tau_i = {\rm Id} - \sigma_i :
E |_{U_i} \>>> S |_{U_i} 
$$
be some splitting of $u$ on $U_i$. Then 
\begin{align}
\tau_j &= \tau_i + u_{ij} \circ \pi \\
\sigma_j &= \sigma_i - \iota \circ u_{ij} \notag 
\end{align}
on $ U_{ij}$. Define $\nabla_i = d_S \circ \tau_i + \sigma_i  \circ d_Q$. Then
$$\nabla_j - \nabla_i \in H^0 (U_{ij} , \Omega^{1}_{X} \otimes
\sH om (Q, S))
$$
is a cocycle. Another choice of $\sigma_i$
verifies 
\begin{align}
\sigma'_i &= \sigma_i - \iota \circ u_i \notag \\
\tau'_i &= \tau_i + u_i \circ \pi \notag
\end{align}
for some $u_i \in H^0 (U_i , \sH om (Q, S))$. Thus 
\begin{gather} \label{well}
\nabla'_i - \nabla_i - d_S (u_i \circ \pi) - \iota \circ u_i \circ
d_Q  \\
= d (u_i ) \in H^0 (U_i , \Omega^{1}_{X} \otimes \sH om (Q, S)),
\notag
\end{gather}
and therefore the class $\alpha_{ij}$ of 
$$
\nabla_j - \nabla_i \in H^1 (X, \Omega^{1}_{X} \otimes \sH om
(Q, S))
$$ is well defined. If this class vanishes, then in a refinement
of $(U_i)$ there are forms $A_i \in H^0 (U_i, \Omega^{1}_{X}
\otimes \sH om (Q, S))$ such that $\nabla_j - \nabla_i = A_i -
A_j$, thus $\nabla = \nabla_i + A_i$ is globally defined and
$\alpha_{ij}$ is the exact obstruction to the existence of
$\nabla$.

On the other hand, the computation in \ref{well}, with $u_i$
replaced by $u_{ij}$, shows at the same time that $\alpha_{ij} =
d u_{ij}$. 
\end{proof} 
Let $X$ be a projective curve over $\C$ and $E$ be a semistable
bundle of degree $0$ on $X$. Then there is a unique 
filtration, which we call {\it the canonical filtration} of $E$,
verifying 
\begin{gather} \label{socle} 
0= E_0 \subset E_1 \subset \ldots \subset E_m = E \\
gr_i E = E_i / E_{i-1} = \ \mbox{socle of} \ E/E_{i-2} .\notag
\end{gather} 
Recall that the {\it socle} of $E$ is the maximal semistable
subbundle of $E$ which splits as a sum $\bigoplus_\nu V_\nu$ of stable
ones. 
\begin{gather} \label{hom} 
\Hom (gr_i E, E/E_i) = \Hom (gr_i E , gr_{i+1} E) \\
= \bigoplus \delta_{\nu\mu} {\rm Id}_{V_\nu} \notag
\end{gather}
with $gr_i E = \bigoplus_\nu V_\nu$, $gr_{i+1} E = \bigoplus_\mu V_\mu$ for
stable bundles $V_\nu$ and $V_\mu$. 

On the other hand, over $\C$, there is a unique unitary
connection $d_i$ on $gr_i E$ by the Narasimhan-Seshadri
correspondence \cite{NS}. 

\begin{prop} \label{exis}
Let $E$ be a semistable bundle of degree 0 on a complex
projective curve, and $E_i$ be its canonical filtration. Then
there is a connection $\nabla$ on $E$ respecting the canonical
filtration and lifting the unitary connections $d_i$ on $E_i/E_{i-1}$. 
\end{prop}

\begin{proof} Since $\sH om (d_m , d_{m-1})$ is unitary, there
is a connection $d_{E/E_{m-1}}$ lifting $d_m$ and $d_{m-1}$ by
lemma \ref{obs}. Assume inductively that $d_{E/E_{\ell}}$ exists.
We want to see that 
$$
d: H^1 (X, \sH om (E/E_{\ell} , gr_{\ell} E)) \>>> H^1 (X, \Omega^{1}_{X}
\otimes \sH om (E/E_{\ell} , gr_{\ell} E)) 
$$
kills the extension of $E/E_{\ell}$ by $gr_{\ell} E$ given by the
canonical filtration, where $d= \sH om (d_{E/E_{\ell}} , d_{\ell} )$. We
show directly that $d$ itself vanishes. Its dual is the
differential 
$$
d^* : H^0 (X, \sH om (gr_{\ell} E, E/E_{\ell})) \>>> H^0 (X, \Omega^{1}_{X}
\otimes \sH om (gr_{\ell} , E/E_{\ell})). 
$$
By the equation \ref{hom}, and the fact that $d^*$ lifts $\sH
om (d_{\ell} , d_{{\ell}+1})$, one has $d^* = \Hom (d_{\ell} ,
d_{{\ell}+1}) =0$. 
\end{proof}

\begin{lem} \label{inv}
Let $X$ be a smooth projective variety defined over a field $k$ of
characteristic zero, $D$ be a smooth irreducible divisor, 
$L$ be an invertible sheaf $L$, and let
$\nabla: L \to \Omega^1_X({\rm log} \ D) \otimes L$ be a
connection. Then the residue
$\res_D(\nabla)$ is $m\cdot {\rm id}$ for a
rational number $m$. Moreover, if $X$ is a curve, $m$ is an integer.
\end{lem}
\begin{proof}
Since $X$ is projective, we may write $L= \sO(A_1-A_2)$ where $A_i$
are smooth divisors meeting transversally. Thus $L$
carries the trivial connection $d_{A}$ with $\res_{A_i}(d_A)=
(-1)^i\cdot{\rm Id}_{L|_{A_I}}$. Hence 
$\omega:=\nabla - d_{A}\in H^0(X, \Omega^1_X({\rm log} \
(A_1+A_2+D))$ with 
$$
m:={\rm res}_D(\omega)= {\rm res}_D \nabla, \ \ {\rm res}_{A_i}
(\omega)= -{\rm res}_{A_i}(d_{A}).
$$
Let $C$ be an ample smooth curve, meeting $D$, $A_1$ and $A_2$
transversally. Then
$$
-(C.A_1)+(C.A_2)+m\cdot (C.D)=
\sum_{q\in C\cap (A_1\cup A_2\cup D)} {\rm res}_{q}(\omega)
=0
$$
and consequently $m \in \Q$ (or $m\in \Z$, if $\dim(X)=1$).
\end{proof}

\begin{lem} \label{sub} 
Let $X$ be a smooth projective variety over a field $k$ of
characteristic zero, let $D = \sum^{\rho}_{i=1} D_i$ be a normal
crossing divisor and 
$$
\nabla : V \to \Omega^{1}_{X} {\rm log} \
D) \otimes V
$$
a connection on a locally free sheaf $V$.
Assume that the eigenvalues of ${\rm res}_{D_i} (\nabla)$ are
zero for $i= 2, \ldots , \rho$ and that the sum of the
eigenvalues of ${\rm res}_{D_1} (\nabla)$ does not lie in $\Q -
\{ 0 \}$ (or not in $\Z -\{0\}$, if $X$ is a curve). Then
$\bigwedge^{\rm max} V$ is numerically trivial. 
\end{lem} 
\begin{proof} 
$\nabla$ induces a connection 
$$
\nabla' : \bigwedge^{\rm max} V \>>> \Omega^{1}_{X} ({\rm log} \
D) \otimes \bigwedge^{\rm max} V.
$$
${\rm res}_{D_i} (\nabla') =0$ for $i = 2, \ldots ,
\rho$, and the image of $\nabla'$ lies in $\Omega^{1}_{X}
(D_1) \otimes \bigwedge^{\rm max} V$. By \ref{inv} ${\rm
res}_{D_1} (\nabla')$ must be a rational number (or an integer),
hence $0$, and $\nabla'$ induces a connection with values in
$\Omega^{1}_{X} \otimes \bigwedge^{\rm max} V$. 
\end{proof} 

\section{Existence of irreducible connections} 

Let $E$ be a semistable bundle of rank $N$ on the curve $X$ and let
$$
\nabla : E \to \Omega^{1}_{X} \otimes E
$$
be a connection. In this section we want to construct a
different connection $\nabla' : E \to \Omega^{1}_{X} ({\rm log}
\ \Sigma) \otimes E$, where $\Sigma = \sum^{\mu}_{i=1} p_i$ is a
reduced divisor in $X$, such that ${\rm Ker} (\nabla' |_{X -
\Sigma})$ is an irreducible local system. If $X$ is defined over
$\C$ this construction and \ref{exis} imply proposition \ref{thm:ev2}.

\begin{prop} \label{twopoints}
Assume that $E$ is not isomorphic to the direct sum $L^{\oplus
N}$ for some $L \in {\rm Pic}^0 (X)$ and let $p,q \in X$ be two
different points. Then there exists $\varphi \in {\rm Hom} (E ,
\Omega^{1}_{X} ({\rm log} (p+q)) \otimes E)$ such that ${\rm
Ker} (\nabla' |_{X- p-q})$ is irreducible for $\nabla' = \nabla +
\varphi$. 
\end{prop}

\begin{proof} 
By assumption there exists a surjection $\tau : E \to S$ for
some bundle $S$ on $X$ of rank $s \geq 2$ such that one of the
following properties holds true: 
\begin{enumerate}
\item[i)] $S$ is stable 
\item[ii)] $S =L_1 \oplus L_2$ for $L_1 \not\cong L_2$ and $L_i
\in {\rm Pic}^0 (X)$ 
\item[iii)] $0 \to T \to S \to L^{\oplus \ell} \to 0$
is an extension of $L^{\oplus \ell}$, for $L \in {\rm Pic}^0
(X)$ with a stable bundle $T$, such that the induced map
$$
H^0 (X, \sO^{\oplus \ell}_{X} ) \>>> H^1 (X, T \otimes L^{-1} ) 
$$
is injective. 
\end{enumerate}
In fact, let $F^{*}_{0} = \{ 0 \} \subset F^{*}_{1} \subset
\ldots \subset F^{*}_{m} = E^*$ be the canonical filtration of
the dual bundle and 
$$
F_0 = \{ 0 \} \subset F_1 = (E^* /F^{*}_{m-1})^* \subset \ldots
\subset F_{m-1} \subset (E^* /F^{*}_{1} )^* \subset F_m = E $$
the dual filtration. If $F_m / F_{m-1}$ contains no semistable
bundle $S$ as in i) or ii) it is a direct sum $L^{\oplus
\ell'}$, for some $\ell' \geq 1$. In this case, 
$$
F_{m-1} / F_{m-2} \>>> E/F_{m-2} \>>> L^{\oplus \ell'} 
$$
is a non-trivial extension and for each direct factor $T$ of
$F_{m-1} /F_{m-2}$ one obtains a surjection from $E$ to a
non-trivial extension 
$$
0 \>>> T \>>> S' \>>> L^{\oplus \ell'}\>>> 0 . 
$$
Leaving out direct factors of $S'$, which are isomorphic to $L$,
one obtains $S$ as in iii). 

For any bundle $F$ on $X$ write $F_q=F\otimes k(q)$.
In order to construct a basis of $E_q$ we fix a basis of $S_q$,
case by case: 
\begin{enumerate}
\item[i)] $\bar{v}_1, \ldots , \bar{v}_{m-1}, \bar{v}_N$ is any
basis of $S_q$. 
\item[ii)] $\bar{v}_1, \bar{v}_N$ is a basis of $S_q$ with
$\bar{v}_1 \not\in (L_i)_q$, for $i = 1,2$. 
\item[iii)] $\bar{v}_1, \ldots , \bar{v}_{m-1} , \bar{v}_N$ is a
basis of $S_q$, such that $T_q \not\subset <\bar{v}_1, \ldots ,
\bar{v}_{m-1} > $. 
\end{enumerate}
Let $K = {\rm Ker} (\tau : E \to S)$ and 
$$
0 \>>> K_q \>>> E_q \> \tau_q >> S_q \>>> 0
$$
the induced sequence of vector spaces. 

Let $v_m, \ldots , v_{N-1}$ be a basis of $K_q$,
and $v_j \in \tau^{-1}_{q} (\bar{v}_j)$, for $j = 1, \ldots , m-1, N$.
Then $v_1, \ldots , v_N$ is a basis of $E_q$. By Serre duality
$$
h^1 (X, \sE nd (E) \otimes \Omega^{1}_{X} ({\rm log} \ p)) = h^0
(X, \sH om (E, E (-p)) = 0 , 
$$
hence the residue map 
$$
H^0 (X, \sE nd (E) \otimes \Omega^{1}_{X} ({\rm log} (p+q))) \>
{\rm res}_q >> {\rm End} (E_q)
$$
is surjective. Choose $\varphi \in {\rm End}
(E, \Omega^{1}_{X} ({\rm log} (p+q)) \otimes E)$ such that ${\rm
res}_q (\varphi)$ is one Jordan block for the eigenvalue $0$,
with respect to $v_1 , \ldots , v_N$ . In particular, the only
${\rm res}_q (\varphi)$ invariant subspaces of $E_q$ are of the form 
${\rm Ker} ({\rm res}_q (\varphi)^{\iota})$. 

Let $\lambda_1, \ldots , \lambda_{\nu}$ be the eigenvalues of
${\rm res}_p (\varphi)$. Replacing $\varphi$ by $\pi \cdot
\varphi$ for some $\pi \not\in \Q (\lambda_1, \ldots ,
\lambda_{\nu})$ we may assume that no linear combination $\Sigma
\rho_i \lambda_i \in \Q - \{ 0 \}$ for $\rho_i \in \Q$. 

Let $V \subset E$ be a subbundle such that $\nabla' (V) \subset
\Omega^{1}_{X} ({\rm log} (p+q)) \otimes V$, for $\nabla' =
\nabla + \varphi$. By \ref{sub} ${\rm deg(V)} =0$, hence $V$ is a
semistable subbundle of $E$, and the image $B$ of
$V$ in $S$ is zero or a semistable subbundle of $S$.

Since ${\rm res}_q (\nabla') = {\rm
res}_q (\varphi)$, for some $\iota \geq 1$
$$
V_q= {\rm Ker} ({\rm res}_q(\varphi)^{\iota}) = < v_1, \ldots ,
v_{\iota} > \subset E_q.
$$ 
In particular $B \neq 0$. Obviously $B=S$ in case i). In case ii) we
remark that $v_1 \in B_q$ and obtain $B=S$, as well. 

If in case iii) $B\neq S$, then $B_q = < v_1, \ldots , v_{\iota} >$ for
$\iota \leq m-1$ and $B \cap T \neq T$. Since the degree of
$B$ is zero, and since $B/(B\cap T) \subset L^{\oplus \ell}$,
$B\cap T = 0$. Then $B \simeq L^
{\oplus \iota}$ and the composite 
$$
H^0 (X, B \otimes L^{-1}) \hookrightarrow H^0 (X, \sO^{\oplus
\ell}_{X} ) \>>> H^1 (X, T \otimes L^{-1})
$$
zero, contradicting the assumptions made. 

Hence $B=S$ in all cases, and $v_n \in V_q$. Therefore $V_q =
E_q$ and $V = E$. 
\end{proof}

If $E = L^{\oplus N}$, then in order to find some
$\varphi$, with ${\rm Ker} (\nabla + \varphi |_{X- \Sigma})$
irreducible, one needs three points $p, q_1, q_2$. In fact,
choosing the ``canonical'' basis $v^{(i)}_{1}, \ldots
v^{(i)}_{N}$ in $E_{q_i}$, induced by the direct sum
decomposition, one has again a surjection 
$$
{\rm End} (E, \Omega^{1}_{X} ({\rm log} (p+q_1 +q_2))) \>>> M (N
\times N , \C) \oplus M (N \times N, \C). 
$$
Let us choose two nilpotent matrices $M_1$ and $M_2$ with
$M^{N-1}_{i} \neq 0$ in such a way, that the (unique)
eigenvector of $M_1$ does not lie in ${\rm Ker} (M^{N-1}_{2})$.
Repeating the argument used in the proof of \ref{twopoints} one
obtains: 

\begin{prop} \label{threepoints} Let $\Sigma = q_1 + q_2 +p$ be
a reduced divisor and $E$ be a semistable bundle
with connection $\nabla$. Then for some $\varphi \in {\rm Hom} (E,
\Omega^{1}_{X} ({\rm log} \ \Sigma ) \otimes E)$ the local
system ${\rm Ker} ((\nabla + \varphi) |_{X - \Sigma} )$ is
irreducible. 
\end{prop}

Under stronger condition on the structure of $E$, it is possible
to choose $\Sigma =p$, as we illustrate in two examples on an
elliptic curve $X$. 

\begin{ex}\label{exone}
Let $L \in {\rm Pic}^0 (X)$, $L \neq \sO$, $E = L \oplus \sO$.
Take $\Sigma = \{ p \}$ a point. Then choose
$$
\nabla = d + \left( \begin{array}{ll} \alpha & \beta \\ \gamma &
\delta \end{array} \right)
$$
where $d $ is the sum of the unitary connections on $L$ and
$\sO$, 
\begin{gather*}
\alpha, \delta \in 
H^0 (X, \Omega_X^1 ({\rm log} \ \Sigma) )=H^0 (X, \Omega_X^1 )  \\
\gamma \in H^0 (X, L^{-1} \otimes \Omega_X^1  ({\rm log} \
\Sigma) )  - H^0 (X, L^{-1} \otimes \Omega_X^1) \\
\beta \in H^0 (X, L \otimes \Omega_X^1  ({\rm log} \
\Sigma) ) -H^0 (X, L^{-1} \otimes \Omega_X^1) .
\end{gather*}
Assume ${\rm res}_q \gamma = \lambda, {\rm res}_q \beta = \mu
$ are chosen such that $x^2 -\lambda \cdot \mu$ has no zero in $\Q$.
If $V\subset E$ of rank 1 is stabilized
by $\nabla$, then ${\rm residue}_p(\nabla|_V)
\not \in \Q$. This contradicts lemma \ref{inv}.
\end{ex}
\begin{ex}\label{extwo}
Let $X$ be an elliptic curve and
$$
0 \>>> \sO_S \> \iota >> E \> \pi >> \sO_Q \>>> 0 
$$
be the non-trivial extension of $\sO_X$ by $\sO_X$.
As we have seen in 3.2, there exists a connection $\nabla$ on
$E$, lifting $d: \sO_X \to \Omega^{1}_{X}$. As
$$
h^0(X,\sE nd(E))=h^1(X,\sE nd(E))=2 \mbox{ \ and \ }
h^1 (X, \sE nd(E (p))) = 0
$$
for any point $p$, whereas $H^0(X,\sO_X)=H^0(X,\sO(p))$, the image of
$$
\res_p: \Hom(E,E(p))=\Hom(E,E\otimes \Omega_X^1(\log \ p))
\>>> M(2\times 2,k)
$$
is a two-dimensional space of matrices of trace $0$.
In particular the image contains some
lower triangular matrix 
$$
M=\left(\begin{array}{cc} \alpha & 0 \\ \gamma & -\alpha \end{array}
\right) \neq 0,
$$
with respect to a basis $v_1,v_2$ with $v_1\in \iota(k(p))$.
Choose $\phi \in \Hom(E,E(p))$ and $\lambda\in k$ with
${\rm res}_p \phi = \lambda\cdot M$, such that $\lambda\alpha \not\in 
\Z - \{0\}$.
By \ref{sub} a rank $1$ subbundle $V \subset E$ with $\nabla (V)
\subset \Omega^{1}_{X} ({\rm log} \ \Sigma ) \otimes V$ is numerically
trivial, hence equal to $\iota(O_X)$. Then $\alpha$ and $\gamma$
are both zero, contradicting the assumption $M\neq0$. 
\end{ex}
\bibliographystyle{plain}

\begin{thebibliography}{99}
\bibitem{Bo} Anosov, D. A.; Bolibruch, A.A.: The Riemann-Hilbert
Problem, Aspects of Mathematics {\bf 22} (1994), Vieweg Verlag. 
\bibitem{G} Gabber, O: letter to A. Beauville, March 1993. 
\bibitem{HN} Harder, G.; Narasimhan, M. S.: On the cohomology
groups of moduli spaces, Math. Ann. {\bf 212} (1975) 215 - 248. 
\bibitem{Le} Lekaus, S.: Diplomarbeit, Universit\"at Essen 1998.
\bibitem{NS} Narasimhan, M. S.; Seshadri, C. S.: Stable and
unitary bundles on a compact Riemann surface, Ann. Math. {\bf
82} (1965), 540 - 567.
\bibitem{Si} Simpson, C.: Higgs bundles and local systems, Publ.
Math. IHES {\bf 75} (1992), 5-95. 
\end{thebibliography}
\renewcommand\refname{References}
 
\end{document}